\documentclass[peerreview, draftclsnofoot, onecolumn]{IEEEtran}

\usepackage{cite} 

\usepackage[cmex10]{amsmath}
\usepackage{amssymb, amsthm, eucal}
\usepackage{array}

\usepackage{graphicx, epsfig}
\usepackage[tight, footnotesize]{subfigure}

\newtheorem{proposition}{Proposition}
\newtheorem{example}{Example}

\def\e{\epsilon}
\def\g{\gamma}
\def\re{\Re}

\def\pttj[#1]{\frac{\partial #1}{\partial \theta_j}}
\def\tpttj[#1]{\tfrac{\partial #1}{\partial \theta_j}}
\def\tpttk[#1]{\tfrac{\partial #1}{\partial \theta_k}}
\def\tpttjk[#1]{\tfrac{\partial^2 #1}{\partial \theta_j \partial \theta_k}}

\def \a {\alpha}
\def\d{\delta}

\def \m {\mu}

\def\argmin{\mathop {\arg \min}}

\def\smskip{\smallskip}
\def\pn{\par\noindent}
\def\lf{\left}
\def\ri{\right}

\def\old#1{}

\def\eqref#1{(\ref{#1})}

\def\jstar{J^*}
\def\p{\pi}
\def\P{\Pi}

\begin{document}

\title{Proper Policies in Infinite-State Stochastic Shortest Path Problems}

\author{Dimitri~P.~Bertsekas%
\thanks{}%
\thanks{D.\ P.\ Bertsekas is with the Computer Information and Decision Science and Engineering Dept., Arizona State University, Tempe, AZ, and the Laboratory for Information and Decision
Systems (LIDS), M.I.T.  Email: dimitrib@mit.edu.
}%
}

\markboth{}%
{Bertsekas: Proper Policies in Infinite-State Stochastic Shortest Path Problems}

 \maketitle

\begin{abstract}
We consider stochastic shortest path  problems with infinite state and control spaces, a nonnegative cost per stage, and a termination state. We extend the notion of a proper policy, a policy that terminates within a finite expected number of steps, from the context of finite state space to the context of infinite state space. We consider the optimal cost function  $\jstar$, and the optimal cost function $\hat J$ over just the proper policies. We show that $\jstar$ and $\hat J$ are the smallest and largest solutions of Bellman's equation, respectively, within a suitable class of Lyapounov-like functions. If the cost per stage is bounded, these functions are those that are bounded over the effective domain of $\hat J$. The standard value iteration algorithm may be attracted to either $\jstar$ or $\hat J$, depending on the initial condition. 
\end{abstract}

\section{Introduction} \label{sec-intro}

In this paper we consider a stochastic discrete-time infinite horizon optimal control problem involving the system
\begin{equation} \label{eq-docsys}
x_{k+1}=f(x_k,u_k,w_k),\qquad k=0,1,\ldots,
\end{equation}
where $x_k$ and $u_k$ are the state and control at stage $k$, which belong to sets $X$ and $U$, $w_k$ is a random disturbance that takes values in a countable set $W$ with given probability distribution $P(w_k\mid x_k,u_k)$, and $f:X\times U\times W\mapsto X$ is a given function. The state and control spaces $X$ and $U$ are arbitrary, but we assume that $W$ is countable to bypass the complicated mathematical measurability issues in the choice of control.
\footnote{The nature of these difficulties is well-documented; see the monograph by Bertsekas and Shreve [1], and the paper by James and Collins [2], which treats stochastic shortest path problems. It may be reasonably conjectured that our analysis can be extended to hold within an appropriate measurability framework, but this undertaking is beyond the scope of the present paper.} 
The control $u_k$ must be chosen from a constraint set $U(x_k)\subset U$ that may depend on the current state $x_k$. The expected cost for the $k$th stage, $E\big\{g(x_k,u_k,w_k)\big\}$, is assumed real-valued and nonnnegative:
\begin{equation} \label{eq-nonnegcost}
0\le E\big\{g(x_k,u_k,w_k)\big\}<\infty,\quad \forall\ x\in X,\ u\in U(x).
\end{equation}
We assume that $X$ contains a special cost-free and absorbing state $t$, referred to as the {\it destination\/}:
\begin{equation} \label{eq-absorbe}
f(t,u,w)=t,\qquad g(t,u,w)=0,\qquad \forall\ u\in U(t),\ w\in W.
\end{equation}
The essence of the problem is to reach or approach the destination with minimum expected cost.

We are interested in   policies of the form $\p=\{\m_0,\m_1,\ldots\}$, where each $\m_k$ is a function mapping $x\in X $ into the control $\m_k(x)\in U(x)$. The set of all  policies is denoted by $\Pi$. Policies of the form $\p=\{\m,\m,\ldots\}$ are called {\it stationary\/}, and will be denoted by $\m$, when confusion cannot arise. 

Given an initial state $x_0$, a policy $\p=\{\m_0,\m_1,\ldots\}$ when applied to the system \eqref{eq-docsys}, generates a random sequence of state-control pairs $\big(x_k,\m_k(x_k)\big)$, $k=0,1,\ldots,$ with cost 
$$J_\p(x_0)= \sum_{k=0}^{\infty}
E^{\p}_{x_0}\Big\{g\big(x_k,\mu_k(x_k),w_k\big)\Big\},\qquad x_0\in X ,$$
where $E^{\p}_{x_0}\{\cdot\}$ denotes expectation with respect to the probability measure corresponding to initial state $x_0$ and policy $\p$, and the series converges in view of the nonnegativity of cost per stage $g$. We view $J_\p$ as a function over $X $, and we refer to it as the cost function of $\p$. For a stationary policy $\m$, the corresponding  cost function is denoted by $J_\m$.
The optimal cost function is defined as
$$\jstar(x)=\inf_{\p\in\P}J_\p(x),\qquad x\in X ,$$
and a policy $\p^*$ is said to be optimal if 
$J_{\p^*}(x)=\jstar(x)$ for all $x\in X .$ We refer to the problem of finding $\jstar$ and an optimal policy as the {\it stochastic shortest path problem} (SSP problem for short). 
We denote by ${\cal E}^+(X)$ the set of functions $J:X\mapsto[0,\infty]$. All equations, inequalities, limit and minimization operations involving functions from this set are meant to be pointwise.  In our analysis, we will use the set of functions
$${\cal J}=\big\{J\in{\cal E}^+(X) \mid J(t)=0\big\}.$$
Since $t$ is cost-free and absorbing, this set contains the cost functions $J_\p$ of all $\p\in \P$, as well as $\jstar$.

It is well known that when $g\ge0$, $\jstar$ satisfies the Bellman equation given by
\begin{equation} \label{eq-bellmaneq}
J(x)=\inf_{u\in U(x)}E\Big\{g(x,u,w)+J\big(f(x,u,w)\big)\Big\},\qquad x\in X,
\end{equation}
where the expected value is with respect to the distribution  $P(w\mid x,u)$. Moreover, an optimal stationary policy (if it exists) may be obtained through the minimization in the right side of this equation when $J=\jstar$ [cf.\ Prop.\ 1(d) in the next section]. One hopes to obtain $\jstar$ in the limit by means of value iteration (VI for short), which starting from some function $J_0\in {\cal J}$, generates a sequence $\{J_k\}\subset  {\cal J}$ according to 
\begin{equation} \label{eq-vieq}
J_{k+1}(x)=\inf_{u\in U(x)}E\Big\{g(x,u,w)+J_k\big(f(x,u,w)\big)\Big\},\quad x\in X.
\end{equation}
However, $\{J_k\}$ may not always converge to $\jstar$ because, among other reasons, Bellman's equation may have multiple solutions within ${\cal J}$.

In a recent paper [3] we have addressed the connections between stability and  the solutions of Bellman's equation in the context of undiscounted discrete-time deterministic optimal control with a termination state. 
In this paper we address similar issues in the context of SSP problems but we focus attention on proper policies, which are the ones that are guaranteed to reach the termination state within a finite expected number of steps, starting from the states where the optimal cost is finite (a precise definition is given in the next section). Proper policies may be viewed as the analog of stable policies in a deterministic context, and their significance is well known in finite-state SSP problems (see e.g., the books [4], [5], [6], [7], [8], [9], [10], and [11], and the references quoted there). For the case where $g\ge0$, the paper by Bertsekas and Tsitsiklis [12] provides an analysis that bears similarity with the one of the present paper, but assumes a finite state space and that there exists an optimal policy that is proper.
In the infinite-state context of this paper and under weaker assumptions, we show that $\hat J$, the optimal cost function over just the proper policies, is the largest solution of Bellman's equation within a set of functions $\widehat {\cal W}\subset {\cal J}$ that majorize $\hat J$, and that the VI algorithm converges to $\hat J$ starting from a function in $\widehat {\cal W}$. Our line of analysis draws its origin from concepts of regularity introduced by the author in the monograph [13] and the paper [14].

To compare our analysis with the existing literature, we note that proper policies for infinite-state SSP problems have  been considered earlier, notably in the works of Pliska [15], and James and Collins [2], where they are called {\it transient\/}. There are a few differences between the frameworks of [15], [2] and this paper, which impact on the results obtained. In particular, the paper [15] uses a similar (but not identical) definition of properness to the one of the present paper, but assumes that all policies are proper, that $g$ is bounded, and that $\jstar$ is \ real-valued. The paper [2] uses the properness definition of [15], and extends the analysis of [11] from finite state space to infinite state space (addressing also measurability issues). Moreover, [2] allows the cost per stage $g$ to take both positive and negative values. However, [2] uses assumptions that guarantee that improper policies cannot be optimal and that $\jstar=\hat J$, while $\jstar$ is real-valued. This is the most important difference from the analysis of this paper. 

Our analysis is also related to the one of Bertsekas and Yu [16], where the case $\jstar\ne\hat J$ was analyzed using perturbation ideas that are similar to the ones of Section 3. The paper [16] assumes that  the state space is finite and that $\jstar$ is real-valued, but allows $g$ to take negative values. Moreover [16] gives an example showing that $\jstar$ may not be a solution of Bellman's equation if improper policies can be optimal, and $g$ can take both positive and negative values. The extension of our results to SSP problems where $g$ takes both positive and negative values may be possible, but our line of analysis relies strongly on the nonnegativity of $g$ and cannot be extended without major modifications. 

\section{Proper Policies and the Perturbed Problem} \label{sec-properpol}

In this section, we will lay the groundwork for our analysis and introduce the notion of a proper policy. To this end, we will use some classical results for stochastic optimal control with nonnegative cost per stage, which stem from the original work of Strauch [17]. For textbook accounts we refer to [1], [8], [11], and for a more abstract development, we refer to the monograph [13].
The following proposition gives the results that we will need. 

\begin{proposition}\label{prp-propnegdp}
The following hold:
\begin{itemize}
\item [(a)]  $\jstar$ is a solution of Bellman's equation and  if $J\in{\cal E}^+(X)$ is another solution, i.e., $J$ satisfies for all $x\in X$
\begin{equation} \label{eq-bellmaneq}
J(x)= \inf_{u\in U(x)}E\Big\{g(x,u,w)+J\big(f(x,u,w)\big)\Big\},
\end{equation}
 then  $\jstar\le J$.
\item [(b)] For all stationary policies $\m$, $J_\m$ is a solution of the equation 
\begin{equation*}
J(x)=E\Big\{g\big(x,\m(x),w\big)+J\big(f\big(x,\m(x),w\big)\big)\Big\},
\end{equation*}
and if $J\in{\cal E}^+(X)$ is another solution, then $J_\m\le J$.
\item [(c)] For every $\e>0$ there exists an $\e$-optimal  policy, i.e., a policy $\p_{\e}$ such that for all $x\in X$, we have
\begin{equation*}
J_{\p_\e}(x)\le \jstar(x)+\e,\qquad\forall\ x\in X.
\end{equation*}
\item [(d)] A stationary policy $\m^*$ is optimal if and only if for all $x\in X$, we have
\begin{equation*}
\m^*(x)\in \argmin_{u\in U(x)}E\Big\{g(x,u,w)+\jstar\big(f(x,u,w)\big)\Big\}.
\end{equation*}
\item [(e)] If $U(x)$ is finite for all $x\in X$, then $J_k\to\jstar$, where $\{J_k\}$ is the sequence generated by the VI algorithm \eqref{eq-vieq}\ starting from any $J_0$ with $0\le J_0\le \jstar$.
\end{itemize}
\end{proposition}

\proof  See [1], Props.\ 5.2, 5.4, and 5.10, or [11], Props.\ 4.1.1, 4.1.3, 4.1.5, 4.1.9. \qed
\vspace{1pc}

For a given state $x\in X$, a policy $\p$ is said to be {\it proper at $x$} if
\begin{equation} \label{eq-propercond}
J_\p(x)<\infty,\qquad \sum_{k=0}^\infty r_k(\p,x)<\infty,
\end{equation}
where $r_k(\p,x_0)$ is the probability that $x_k\ne t$ when using $\p$ and starting from $x_0=x$. Note that the sum $\sum_{k=0}^\infty r_k(\p,x)$ is the expected number of steps to reach the destination starting from $x$ and using $\p$. 

We denote by  $\widehat{\P}_x$ the set of all policies that are proper at $x$, and we use the notation
\begin{equation} \label{eq-defcalc}
{\cal C}=\big\{(\p,x)\mid \p\in \hat \P_x\big\}.
\end{equation}
We denote by $\hat J$ the corresponding restricted optimal cost function,
\begin{equation*}
\hat J(x)=\inf_{(\p,x)\in {\cal C}} J_\p(x)=\inf_{\p\in \widehat{\P}_x} J_\p(x),\qquad x\in X.
\end{equation*}
Finally we denote by $\widehat{X}$ the effective domain of $\hat J$, i.e.,
\begin{equation} \label{eq-effdomain}
\widehat{X}=\big\{x\in X\mid \hat J(x)<\infty\big\}.
\end{equation}
Note that $\widehat{X}$ is the set of all $x$ such that $\widehat \P_x$ is nonempty.

The definition of proper policy just given differs from the definition of a transient policy adopted by James and Collins [2]. In particular, the definition of [2] requires that the expected number of steps to reach the destination is uniformly bounded over the initial state $x$ (see  [2], p.\ 608) and is not tied to a single state $x$.

For any $\d>0$, let us consider the $\d$-{\it perturbed optimal control problem\/}. This is the same problem as the original, except that the cost per stage is changed to
\begin{equation*}
g(x,u,w)+\d,\qquad \forall\ x\ne t,
\end{equation*}
while $g(x,u,w)$ is left unchanged at 0 when $x=t$. Thus $t$ is still cost-free as well as absorbing in the $\d$-perturbed problem. The $\d$-perturbed cost function of a policy $\p$ is denoted by $J_{\p,\d}$ and is given by
\begin{equation} \label{eq-costexpr}
J_{\p,\d}(x)=J_\p(x)+ \d\sum_{k=0}^\infty r_k(\p,x).
\end{equation}
We denote by $\hat J_\d$ the optimal cost function of the $\d$-perturbed problem, i.e.,
$\hat J_\d(x)=\inf_{\p\in\Pi}J_{\p,\d}(x)$. The following proposition relates the $\d$-perturbed problem with proper policies.

\begin{proposition}\label{prp-propdeltaper}
\begin{itemize}
\item [(a)]  A policy $\p$ is proper at a state $x\in X$ if and only if $J_{\p,\d}(x)<\infty$.
\item [(b)] We have  $\hat J_{\d}(x)<\infty$ for all $\d>0$ if and only if $x\in \widehat{X}$.
\item [(c)] For every $\e>0$ and $\d>0$, there exists a policy $\p_\e$ that is proper at all $x\in \widehat{X}$ and is $\e$-optimal for the $\d$-perturbed problem, i.e.,
\begin{equation*}
J_{\p_\e,\d}(x)\le \hat J_\d(x)+\e,\qquad \forall\ x\in X.
\end{equation*}
\end{itemize}
\end{proposition}

\proof (a) Follows from Eq.\ \eqref{eq-costexpr} and the definition \eqref{eq-propercond} of a proper policy.
\smskip
\pn (b) If $x\in \widehat{X}$ there exists a policy $\p$ that is proper at $x$, and by part (a), $\hat J_\d(x)\le J_{\p,\d}(x)<\infty$ for all $\d>0$. Conversely, if $\hat J_\d(x)<\infty$, there exists $\p$ such that $J_{\p,\d}(x)<\infty$, implying [by part (a)] that $\p\in \widehat{\P}_x$, so that $x\in \widehat{X}$.
\smskip
\pn (c) By Prop.\ \ref{prp-propnegdp}(c), there exists an $\e$-optimal policy $\p_\e$ for the $\d$-perturbed problem, so we have $J_{\p_\e,\d}(x)\le \hat J_\d(x)+\e$ for all $x\in X$. Hence $J_{\p_\e,\d}(x)<\infty$ for all $x\in \widehat{X}$, implying by part (a) that $\p_\e$ is proper  at all $x\in \widehat{X}$.
\qed
\vspace{1pc}

The next proposition shows that the cost function $\hat J_\d$ of the $\d$-perturbed problem can be used to approximate $\hat J$. 

\begin{proposition}\label{prp-proptaffmonf}
We have 
$\lim_{\d\downarrow0}\hat J_\d(x)=\hat J(x)$ for all $x\in X.$
\end{proposition}

\proof Let us fix $\d>0$, and  for a given $\e>0$, let $\p_\e$ be a policy that is proper at all $x\in \widehat{X}$ and is $\e$-optimal for the $\d$-perturbed problem [cf.\ Prop.\ \ref{prp-propdeltaper}(c)]. By using Eq.\ \eqref{eq-costexpr}, we have  for all $\e>0$, $x\in \widehat{X}$, and $\p\in \widehat{\P}_x$,
\begin{align*}
\hat J(x)-\e&\le J_{\p_\e}(x)-\e\\
&\le J_{\p_\e,\d}(x)-\e\\
&\le  \hat J_\d(x)\\
&\le J_{\p,\d}(x)\\
&= J_{\p}(x)+w_{\p,\d}(x),\qquad \forall\ x\in \widehat{X},
\end{align*}
where 
$$w_{\p,\d}(x)=\d\sum_{k=0}^\infty r_k(\p,x)<\infty,\qquad \forall\ x\in \widehat{X}.$$
By taking the limit as $\e\downarrow 0$, we obtain for all $\d>0$ and $\p\in \widehat{\P}_x$,
$$\hat J(x)\le \hat J_\d(x)\le J_{\p}(x)+w_{\p,\d}(x),\qquad \forall\ x\in \widehat{\P}_x.$$
We have $\lim_{\d\downarrow0}w_{\p,\d}(x)=0$ for all $x\in \widehat{X}$ and $\p\in \widehat{\P}_x$, so  by taking the limit as $\d\downarrow0$ and then the infimum over all  $\p\in \widehat{\P}_x$,
$$\hat J(x)\le \lim_{\d\downarrow0}\hat J_\d(x)\le \inf_{\p\in \widehat{\P}_x}J_\p(x)=\hat J(x),\qquad \forall\ x\in \widehat{X},$$
from which $\hat J(x)= \lim_{\d\downarrow0}\hat J_\d(x)$  for all $x\in \widehat{X}$. Since by Prop.\ \ref{prp-propdeltaper}(b), we also have $\hat J_\d(x)=\hat J(x)=\infty$ for all $x\notin \widehat{X}$, the result follows. 
\qed
\vspace{1pc}

\section{Main Results}

By Prop.\ \ref{prp-propnegdp}(a), $\hat J_\d$ solves Bellman's equation for the $\d$-perturbed problem, while by Prop.\ \ref{prp-proptaffmonf}, $\lim_{\d\downarrow0}\hat J_\d(x)=\hat J(x)$. This suggests that $\hat J$ solves the unperturbed Bellman equation, which is the ``limit" as $\d\downarrow0$ of the $\d$-perturbed version. Indeed,  under a certain assumption we will show a stronger result, namely that $\hat J$ is the unique solution of Bellman's equation within the set of functions
\begin{equation} \label{eq-setcalw}
\widehat{\cal W}=\Big\{J\in {\cal J}\mid \hat J\le J,\,E^\p_{x_0}\big\{J(x_k)\big\}\to0,\ \forall\ (\p,x_0)\in{\cal C}\Big\},
\end{equation}
where 
${\cal C}$ is given by
Eq.\ \eqref{eq-defcalc}, $E^\p_{x_0}\{\cdot\}$ denotes expected value with respect to the probability measure corresponding to initial state $x_0$ under policy $\p$, and $E^\p_{x_0}\big\{J(x_k)\big\}$ denotes the expected value of the function $J$ along the sequence $\{x_k\}$ generated starting from $x_0$ and using $\p$. The functions in $\widehat{\cal W}$ are the ones whose expected value is decreasing to 0 along the trajectories generated by the proper policies, so they may be interpreted as a type of Lyapounov functions.

Given a policy $\p=\{\m_0,\m_1,\ldots\}$,  we denote by $\p_k$ the policy 
\begin{equation} \label{eq-tailpol}
\p_k=\{\m_k,\m_{k+1},\ldots\}.
\end{equation}
We first show a preliminary result. 

\begin{proposition}\label{prp-propprelim}
The following hold:
\begin{itemize}
\item [(a)]   For all pairs $(\p,x_0)\in {\cal C}$ and $k=0,1,\ldots$, we have
$$0\le E_{x_0}^{\p}\big\{\hat J(x_k)\big\}\le E_{x_0}^{\p}\big\{J_{\p_k}(x_k)\big\},$$
where $\p_k$ is the policy given by Eq.\ \eqref{eq-tailpol}.
\item [(b)] The set $\widehat{\cal W}$ of Eq.\ \eqref{eq-setcalw} contains $\hat J$, as well as all  $J\in \widehat{\cal W}$ satisfying 
$\hat J\le J\le c\hat J$ 
for some $c\ge1$.
\end{itemize}
\end{proposition}

\proof (a) For any pair  $(\p,x_0)\in {\cal C}$ and $\d>0$, we have
\begin{align*}
J_{\p,\d}(x_0)=E_{x_0}^{\p}\bigg\{J_{\p_k,\d}(x_k)&+k\d\\
&+\sum_{m=0}^{k-1}g\big(x_m,\m_m(x_m),w_m\big)\bigg\}.
\end{align*}
Since $J_{\p,\d}(x_0)<\infty$ [cf.\ Prop.\ \ref{prp-propdeltaper}(a)], it follows that $E_{x_0}^{\p}\big\{J_{\p_k,\d}(x_k)\big\}<\infty$. Hence for all $x_k$ that can be reached with positive probability using $\p$ and starting from $x_0$, we have 
$J_{\p_k,\d}(x_k)<\infty$, implying [by Prop.\ \ref{prp-propdeltaper}(a)] that $(\p_k,x_k)\in {\cal C}$ and hence $\hat J(x_k)\le J_{\p_k}(x_k)$. By applying $E_{x_0}^{\p}\{\cdot\}$ to this last inequality, the result follows. 

\smskip
\pn (b)  
We have for all $(\p,x_0)\in{\cal C}$,
\begin{equation} \label{eq-fjpexpression}
J_\p(x_0)=E_{x_0}^{\p}\Big\{g\big(x_0,\m_0(x_0),w_0\big)\Big\}+E_{x_0}^{\p}\big\{J_{\p_1}(x_{1})\big\},
\end{equation}
 and  for all $m=1,2,\ldots,$
\begin{align} \label{eq-finitejpi}
E_{x_0}^{\p}\big\{J_{\p_m}(x_m)\big\}=&E_{x_0}^{\p}\Big\{g\big(x_m,\m_m(x_m),w_m\big)\Big\}\notag\\
&\ \ \ \ \ \ +E_{x_0}^{\p}\big\{J_{\p_{m+1}}(x_{m+1})\big\},
\end{align}
where $\{x_m\}$ is the sequence generated starting from $x_0$ and using $\p$. By using repeatedly the expression \eqref{eq-finitejpi} for $m=1,\ldots,k-1$, and combining it with Eq.\ \eqref{eq-fjpexpression}, we obtain for all $k=1,2,\ldots,$ and $(\p,x_0)\in {\cal C}$,
$$J_\p(x_0)=E_{x_0}^{\p}\big\{J_{\p_k}(x_k)\big\}+\sum_{m=0}^{k-1}E_{x_0}^{\p}\Big\{g\big(x_m,\m_m(x_m),w_m\big)\Big\}.$$
The rightmost term above tends to $J_\p(x_0)$ as $k\to\infty$, so by using the fact $J_\p(x_0)<\infty$, we obtain 
$$E_{x_0}^{\p}\big\{J_{\p_k}(x_k)\big\}\to0,\qquad \forall\ (\p,x_0)\in {\cal C}.$$
By part (a), it follows that
$E_{x_0}^{\p}\big\{\hat J(x_k)\big\}\to0$ for all $(\p,x_0)\in {\cal C},$
so that $\hat J\in \widehat{\cal W}$. This also implies that 
$E_{x_0}^{\p}\big\{J(x_k)\big\}\to0$ for all $(\p,x_0)\in {\cal C},$
if $\hat J\le J\le c\hat J$ for some $c\ge1$.
\qed
\vspace{1pc}

We can now prove our main result. We denote by $X^*$ the effective domain of $\jstar$:
\begin{equation} \label{eq-effdomain}
{X^*}=\big\{x\in X\mid \jstar(x)<\infty\big\}.
\end{equation}
We also denote by $ Q^*_\d(x,u)$ the optimal Q-factor of a pair $(x,u)$ in the $\d$-perturbed problem:
$$Q^*_\d(x,u)=E\Big\{g(x,u,w)+\d+\hat J_{\d}\big(f(x,u,w)\big)\Big\}.$$

The following proposition shows that $\hat J$ is the unique solution of the Bellman equation within the set $\widehat{\cal W}$ of Lyapounov functions under a certain assumption relating to the states in $X^*$. This assumption can often be easily verified in practice. It is satisfied for example if there exists a policy $\p$ (necessarily proper at all $x\in X^*$) such that $J_{\p,\d}$ is bounded over the set $X^*$. Later, we will also prove the result under the alternative assumption that the set of disturbances $W$ is finite.

\begin{proposition}\label{prp-proptaffmonft}Assume that there exists a $\d>0$ such that
\begin{equation} \label{eq-specialassumption}
Q^*_\d(x,u)<\infty,\qquad \forall\ x\in X^*,\ u\in U(x).
\end{equation}
Then:
\begin{itemize}
\item [(a)]  $\hat J$ is the unique solution of the Bellman Eq.\ \eqref{eq-bellmaneq} within the set $\widehat{\cal W}$ of Eq.\ \eqref{eq-setcalw}.
\item [(b)]  ({\it VI Convergence\/})  If $\{J_k\}$ is the sequence generated by the VI algorithm \eqref{eq-vieq} starting with some 
 $J_0\in \widehat{\cal W}$, then $J_k\to\hat J$.
\item[(c)]  ({\it Optimality Condition\/})  If $\m$ is a stationary policy that is proper at all $x\in \widehat{X}$, and for all $x\in X$ we have
\begin{equation} \label{eq-optcond}
\m(x)\in \argmin_{u\in U(x)}E\Big\{g(x,u,w)+\hat J\big(f(x,u,w)\big)\Big\},
\end{equation}
then $\m$ is optimal over the set of proper policies, i.e., $J_{\m}=\hat J$. Conversely, if $\m$ is optimal within the set of proper policies, then it satisfies the preceding condition \eqref{eq-optcond}.
\end{itemize}
\end{proposition}

\proof (a), (b) By Prop.\ \ref{prp-propprelim}(b), $\hat J\in \widehat{\cal W}$.
We will first show that $\hat J$ is a solution of Bellman's equation and then show that it is the unique solution within $\widehat{\cal W}$ by showing the convergence of VI [cf.\ part (b)].
Since $\hat J_\d$ solves the Bellman equation for the $\d$-perturbed problem, and $\hat J_\d\ge \hat J$ (cf.\ Prop.\ \ref{prp-proptaffmonf}), we have for all $\d>0$ and $x\ne t$,
\begin{align*}
\hat J_\d(x)&=\inf_{u\in U(x)}E\Big\{g(x,u,w)+\d+ \hat J_\d\big(f(x,u,w)\big)\Big\}\\
&\ge \inf_{u\in U(x)}E\Big\{g(x,u,w)+ \hat J_\d\big(f(x,u,w)\big)\Big\}\\
&\ge \inf_{u\in U(x)}E\Big\{g(x,u,w)+ \hat J\big(f(x,u,w)\big)\Big\}.
\end{align*}
By taking the limit as $\d\downarrow0$ and using Prop.\ \ref{prp-proptaffmonf}, we obtain
\begin{equation} \label{eq-forineq}
\hat J(x)\ge \inf_{u\in U(x)}E\Big\{g(x,u,w)+\hat J\big(f(x,u,w)\big)\Big\},\qquad \forall\ x\in X.
\end{equation}

To prove the reverse inequality, we consider two cases:

(1) $x\notin X^*$, i.e., $\jstar(x)=\infty$. Then from Bellman's equation, we have
$$\infty=\jstar(x)=\inf_{u\in U(x)}E\Big\{g(x,u,w)+\jstar\big(f(x,u,w)\big)\Big\}.$$
Since $\hat J\ge \jstar$ it then follows that
\begin{equation} \label{eq-caseinfinity}
\infty=\hat J(x)=\inf_{u\in U(x)}E\Big\{g(x,u,w)+\hat J\big(f(x,u,w)\big)\Big\},\ \forall\ x\notin X^*.
\end{equation}

(2) $x\in X^*$. Then we let $\{\d_m\}$ be a sequence with $\d_m\downarrow0$.  We have for all $m$, $x\ne t$, and $u\in U(x)$,
\begin{align*}
Q^*_{\d_m}(x,u)&=E\Big\{g(x,u,w)+\d_m+\hat J_{\d_m}\big(f(x,u,w)\big)\Big\} \\
&\ge\inf_{v\in U(x)}E\Big\{g(x,v,w)+\d_m+\hat J_{\d_m}\big(f(x,v,w)\big)\Big\}\\
&=\hat J_{\d_m}(x).
\end{align*}
We now take limit as $m\to\infty$ in the preceding relation. The condition \eqref{eq-specialassumption} implies that for all $m$ sufficiently large the left side is finite for each $u\in U(x)$, so we can apply the monotone convergence theorem to interchange limit as $m\to\infty$ and expectation.\footnote{We are using here the following version of the monotone convergence theorem: Let $\{h_m\}$ be a sequence of monotonically nonincreasing functions $h_m:\{1,2,\ldots\}\mapsto\re$, let $\{p_1,p_2,\ldots\}$ be a probability distribution, and assume that for some function $\bar h:\{1,2,\ldots\}\mapsto\re$ such that $h_m(i)\le \bar h(i)$ for all $m$ and $i$, we have
$\sum_{i=1}^\infty p_i\bar h(i)<\infty.$
Then
$$\lim_{m\to\infty}\sum_{i=1}^\infty p_ih_m(i)= \sum_{i=1}^\infty p_i\lim_{m\to\infty} h_m(i).$$
We give the proof, which is simple in the discrete distribution case considered here: Let $h$ be the pointwise limit of $\{h_m\}$, i.e., $h(i)=\lim_{m\to\infty}h_m(i)$ for all $i$. Since $\{h_m\}$ is nonincreasing, we have
$$\sum_{i=1}^\infty p_ih_m(i)\ge \sum_{i=1}^\infty p_ih(i),\qquad \forall\ m=0,1,\ldots,$$
so that
$$\lim_{m\to\infty}\sum_{i=1}^\infty p_ih_m(i)\ge \sum_{i=1}^\infty p_ih(i).$$
Conversely, since $\bar h(i)-h_m(i)\ge0$ for all $m$, we have for every $N\ge 1$
$$\sum_{i=1}^\infty p_i\big(\bar h(i)-h_m(i)\big)\ge \sum_{i=1}^N p_i\big(\bar h(i)-h_m(i)\big),$$
and hence
$$\lim_{m\to\infty}\sum_{i=1}^\infty p_i\big(\bar h(i)-h_m(i)\big)\ge \lim_{m\to\infty}\sum_{i=1}^N p_i\big(\bar h(i)-h_m(i)\big),$$
so that
$$\sum_{i=1}^\infty p_i\bar h(i)-\lim_{m\to\infty}\sum_{i=1}^\infty p_i h_m(i)\ge \sum_{i=1}^N p_i \bar h(i)-\sum_{i=1}^N p_i h(i).$$
By taking the limit as $N\to\infty$ and using the fact that $\sum_{i=1}^\infty p_i\bar h(i)$ is finite (so we can cancel it from both sides of the inequality), we obtain
$$\lim_{m\to\infty}\sum_{i=1}^\infty p_i h_m(i)\le \sum_{i=1}^\infty p_i h(i),$$
thus completing the proof.
} 
Since $\lim_{\d_m\downarrow0}\hat J_{\d_m}=\hat J$ (cf.\ Prop.\ \ref{prp-proptaffmonf}), we obtain
$$E\Big\{g(x,u,w)+\hat J\big(f(x,u,w)\big)\Big\}\ge \hat J(x),\quad \forall\ x\in X^*,\ u\in U(x),$$
so that 
\begin{equation} \label{eq-backineq}
\inf_{u\in U(x)}E\Big\{g(x,u,w)+\hat J\big(f(x,u,w)\big)\Big\}\ge \hat J(x),\qquad \forall\ x\in X^*.
\end{equation}
Thus by combining Eqs.\ \eqref{eq-forineq}, \eqref{eq-caseinfinity}, and \eqref{eq-backineq}, we see that 
$$\hat J(x)=\inf_{u\in U(x)}E\Big\{g(x,u,w)+\hat J\big(f(x,u,w)\big)\Big\},\qquad \forall\ x\in X,$$
and that $\hat J$ is a solution of Bellman's equation.

We will next show that $J_k\to\hat J$ starting from every initial $J_0\in \widehat{\cal W}$ [cf.\ part (b)]. Indeed, for $x_0\in \widehat{X}$ and any  $\p=\{\m_0,\m_1,\ldots\}\in\widehat{\P}_{x_0}$, let $\{x_k\}$ be the generated sequence starting from $x_0$. Since from the definition of the VI sequence $\{J_k\}$ [cf.\ Eq.\ \eqref{eq-vieq}], we have for all $x\in X$, $u\in U(x)$, $k=1,2,\ldots$,
$$J_k(x)\le E\Big\{g(x,u,w)+J_{k-1}\big(f(x,u,w)\big)\Big\},$$
it follows that
$$J_k(x_0)\le E_{x_0}^\p\lf\{J_0(x_k)+\sum_{m=0}^{k-1} g\big(x_m,\m_m(x_m),w_m\big)\ri\}.$$
Since $J_0\in \widehat{\cal W}$, we have $E_{x_0}^\p\big\{J_0(x_k)\big\}\to0$, so by taking the limit as $k\to\infty$ in the preceding relation, it follows that $\limsup_{k\to\infty}J_k(x_0)\le J_\p(x_0)$. By taking the infimum over all $\p\in\widehat{\P}_{x_0}$, we obtain 
$\limsup_{k\to\infty}J_k(x_0)\le \hat J(x_0).$ Conversely, since $\hat J\le J_0$ and $\hat J$ is a solution of Bellman's equation (as shown earlier), it follows by induction that $\hat J\le J_k$ for all $k$. Thus 
$\hat J(x_0)\le \liminf_{k\to\infty}J_k(x_0),$
 implying that $J_k(x_0)\to\hat J(x_0)$ for all $x_0\in \widehat{X}$. We also have $\hat J\le J_k$ for all $k$, so that $\hat J(x_0)=J_k(x_0)=\infty$ for all $x_0\notin \widehat{X}$. This completes the proof of part (b).
Finally, since $\hat J\in \widehat{\cal W}$ and $\hat J$ is a solution of Bellman's equation, part (b) implies the uniqueness assertion of part (a).

\smskip
\pn (c) If $\m$ is proper at all $x\in \widehat{X}$ and Eq.\ 
\eqref{eq-optcond} holds, then 
$$\hat J(x)=E\Big\{g\big(x,\m(x),w\big)+\hat J\big(f(x,\m(x),w)\big)\Big\},\qquad x\in X.$$
By Prop.\ \ref{prp-propnegdp}(b), this implies that  
$J_\m\le \hat J$, so 
$\m$ is optimal over the set of proper policies. Conversely, assume that $\m$ is proper at all $x\in \widehat{X}$ and $J_\m=\hat J$. Then by Prop.\ \ref{prp-propnegdp}(b), we have
$$\hat J(x)=E\Big\{g\big(x,\m(x),w\big)+\hat J\big(f(x,\m(x),w)\big)\Big\},\qquad x\in X,$$
and since [by part (b)] $\hat J$ is a solution of Bellman's equation,
$$\hat J(x)=\inf_{u\in U(x)}E\Big\{g(x,u,w)+\hat J\big(f(x,u,w)\big)\Big\},\qquad x\in X.$$
Combining the last two relations, we obtain Eq.\ \eqref{eq-optcond}.  \qed
\vspace{1pc}

Let us also state our main result under an alternative assumption, which makes the connection with our earlier deterministic results of the paper [3], where an assumption such as Eq.\ \eqref{eq-specialassumption} is not needed.

\begin{proposition}\label{prp-proptaffmonfth}Assume that the disturbance set $W$ is finite. Then the conclusions of Prop.\ \ref{prp-proptaffmonft} hold.
\end{proposition}

\proof The monotone convergence argument for the proof of Eq.\ \eqref{eq-backineq} goes through using the finiteness of $W$ in place of the assumption \eqref{eq-specialassumption}. \qed
\vspace{1pc}

We note that some additional assumption, like Eq.\ \eqref{eq-specialassumption} or the finiteness of $W$, is necessary to prove our results for SSP problems. In this respect, we note that the original version of the proposition, which appeared in the IEEE Trans.\ on Aut.\ Control, was flawed in that it was valid only for the case where $W$ is finite. This was pointed out to us  by Yi Zhang (private communication), who constructed the following example. 

\begin{example}\label{ex-example1}Let $X=\{t, 0,1,2,\ldots\}$, where $t$ is the termination state, and let $g(x,u,w)\equiv0$, so that $\jstar(x)\equiv0$. There is only one control at each state, and hence only one policy. The transitions are as follows:

From each state ${x=2,3,\ldots}$ we move deterministically to state $x-1$, from state 1 we move deterministically to state $t$, and from state 0 we move to state $x=1,2,\ldots$, with probability $p_x$ such that $\sum_{x=1}^\infty x p_x=\infty$ [so at state 0, the assumption \eqref{eq-specialassumption} and the finiteness of $W$ are violated].

Here the unique policy is proper at all $x=1,2,\ldots$, and we have $\skew6\hat J(x)=J^*(x)=0$.  However, the policy is not proper at $x=0$, since the expected number of transitions from $x=0$ to termination is $\sum_{x=1}^\infty x p_x=\infty$. As a result the set $\widehat{\P}_0$ is empty and we have $\hat J(0)=\infty$. Thus $\hat J$ does not satisfy the Bellman equation for $x=0$, since
$$\infty=\hat J(0)\ne E\Big\{g(0,u,w)+\hat J\big(f(0,u,w)\big)\Big\}=\sum_{x=1}^\infty p_x\hat J(x)=0.$$
\end{example}
\smskip

Suppose now that the set of proper policies is sufficient in the sense that it can achieve the same optimal cost as the set of all policies, i.e., $\hat J=\jstar$.  Then, Prop.\ \ref{prp-proptaffmonft} or Prop.\ \ref{prp-proptaffmonfth} (under the corresponding assumptions) imply that $\jstar$ is the unique solution of Bellman's equation within $\widehat{\cal W}$, and the VI algorithm converges to $\jstar$ starting from any $J_0\in \widehat{\cal W}$. Under additional conditions, such as finiteness of $U(x)$ for all $x\in X$ [cf.\ Prop.\ \ref{prp-propnegdp}(e)], the VI algorithm converges to $\jstar$ starting from any $J_0\in{\cal J}$ with $E^\p_{x_0}\big\{J(x_k)\big\}\to0$, for all $(\p,x_0)\in{\cal C}$.

\section{The Multiplicity of Solutions of Bellman's Equation}

Let us now discuss the issue of multiplicity of solutions of Bellman's equation within the set of functions
$${\cal J}=\big\{J\in{\cal E}^+(X) \mid J(t)=0\big\}.$$
We know from Prop.\ \ref{prp-propnegdp}(a) and Prop.\ \ref{prp-proptaffmonft}(a) (or Prop.\ \ref{prp-proptaffmonfth}) that $\jstar$ and $\hat J$ are solutions, and that all other solutions $J$ must satisfy either $\jstar\le J\le \hat J$ or $J\notin\widehat{\cal W}$.

In the special case of a deterministic problem (one where the disturbance $w_k$  takes a single value), it was shown in the paper [3] that $\hat J$ is the largest solution of Bellman's equation within ${\cal J}$, so all solutions $J\in {\cal J}$  satisfy $\jstar\le J\le \hat J$. Moreover, it was shown through examples that there can be  any number of solutions that lie between $\jstar$ and $\hat J$: a finite number, an infinite number, or none at all. 

In stochastic problems, however, the situation is strikingly different. There can be an infinite number of solutions $J\in {\cal J}$ such that $J\ne \hat J$ and $J\ge \hat J$, even when the set $W$ is finite, as illustrated by the following example. Of course, by Prop.\ \ref{prp-proptaffmonft}(a) or Prop.\ \ref{prp-proptaffmonfth}, under the corresponding assumptions, these solutions must lie outside $\widehat{\cal W}$. 

\begin{example}\label{ex-example2}Let $X=\re$, $t=0$, and assume that there is only one control at each state. The disturbance $w_k$ takes two values: 1 and 0 with probabilities $\a\in(0,1)$ and $1-\a$, respectively. The system equation is
$$x_{k+1}={w_k x_k\over \a},$$
and there is no cost at each state and stage [$g(x,u,w)\equiv0$].
Thus from state $x_k$ we move to $x_k/\a$ with probability $\a$ and to the termination state $t=0$ with probability $1-\a$. Here, the only admissible policy is proper, and we have
$$J^*(x)=\skew6\hat J(x)=0,\qquad \forall\ x\in X.$$
Bellman's equation has the form
$$J(x)=(1-\a)J(0)+\a J\left({x\over a}\right),\qquad x\in X,$$
and has an infinite number of solutions within ${\cal J}$ in addition to $J^*$ and $\skew6\hat J$: any positively homogeneous function, such as, for example, $J(x)=\g|x|$, $\g>0$, is a solution. Consistently with Prop.\ \ref{prp-proptaffmonft}(a), none of these solutions belongs to $\widehat{\cal W}$, since $x_k$ is either equal to $x_0/\a^k$ (with probability $\a^k$) or equal to 0  (with probability $1-\a^k$), and, for example, $E\big\{\g |x_k|\big\}=\g |x_0|$ for all $k$.
\end{example}

Let us also note that in the case of linear-quadratic problems, the number of solutions of the Riccati equation has been the subject of considerable investigation, starting with the papers by Willems [18] and Kucera [19], [20], which were followed up by several other papers. These works adopt various assumptions relating to controllability and observability. Because of these assumptions and also because solutions of the Riccati equation give rise to solutions of the Bellman equation, but not reversely, it appears that the full characterization of the set of solutions of the Bellman equation remains an interesting open research question at present.
 
\section{The Case of Bounded Cost per Stage}

Let us consider the special case where the cost per stage $g$ is bounded over $X\times U\times W$, i.e.,
\begin{equation} \label{eq-costbounded}
\sup_{(x,u,w)\in X\times U\times W}g(x,u,w)<\infty.
\end{equation}
We will show that $\hat J$ is the largest solution of Bellman's equation within the class of functions that are bounded over the effective domain $\widehat{X}$ of $\hat J$ [cf.\ Eq.\ \eqref{eq-effdomain}].

We say that a policy $\p$ is {\it uniformly proper} if there is a uniform bound on the expected number of steps to reach the destination from states $x\in\widehat{X}$ using $\p$:
$$\sup_{x\in \widehat{X}}\sum_{k=0}^\infty r_k(\p,x)<\infty.$$
Since we have for all $\p\in\widehat \P_{x_0}$,
$$J_\p(x_0)\le \lf(\sup_{(x,u,w)\in X\times U\times W}g(x,u,w)\ri)\cdot \sum_{k=0}^\infty r_k(\p,x_0)<\infty,$$
it follows that the cost function $J_\p$ of a uniformly proper $\p$ belongs to the set ${\cal B}$, defined by
\begin{equation} \label{eq-setcalb}
{\cal B}=\lf\{J\in {\cal J}\ \Big|\  \sup_{x\in\widehat{X}}J(x)<\infty\ri\}.
\end{equation}
When $\widehat{X}=X$, the notion of a uniformly proper policy coincides with the notion of a transient policy used in [2] and [15], which itself descends from earlier works. However, our definition is somewhat more general, since it also applies to the case where $\widehat{X}$ is a strict subset of $X$.

Let us denote by $\widehat {\cal W}_b$ the set of functions
$$\widehat {\cal W}_b=\{J\in {\cal B}\mid \hat J\le J\}.$$
The following proposition provides conditions for $\hat J$ to be the largest fixed solution of the Bellman equation within ${\cal B}$. Its assumptions include the existence of a uniformly proper policy, which implies that $\hat J$ belongs to ${\cal B}$.

\begin{proposition}\label{prp-1}
Assume that the cost per stage $g$ is bounded over $X\times U\times W$ [cf.\ Eq.\ \eqref{eq-costbounded}], and that there exists a uniformly proper policy. Assume further that Eq.\ \eqref{eq-specialassumption} holds or that the set $W$ is finite. Then:
\begin{itemize}
\item [(a)]  $\hat J$ is the unique solution of the Bellman Eq.\ \eqref{eq-bellmaneq} within the set $\widehat {\cal W}_b$. Moreover, if $\hat J=\jstar$, then $\jstar$ is the unique solution of Bellman's equation within ${\cal B}$.
\item [(b)] If $\{J_k\}$ is the sequence generated by the VI algorithm \eqref{eq-vieq} starting with some 
 $J_0\in {\cal B}$ with $J_0\ge \hat J$, then $J_k\to\hat J$.
 \end{itemize}
 \end{proposition}
 
 \def\tl{\tilde}

\proof Since, as noted earlier, the cost function of a uniformly proper policy belongs to ${\cal B}$, it follows that $\hat J$ also belongs to ${\cal B}$. On the other hand, for all $J\in {\cal B}$, we have
$$E^\p_{x_0}\big\{J(x_k)\big\}\le \lf(\sup_{x\in \widehat{X}}J(x)\ri)\cdot r_k(\p,x_0)\to0,\qquad  \forall\ \p\in\widehat \P_{x_0}.$$
It follows that the set
$\widehat {\cal W}_b$
is contained in $\widehat {\cal W}$, while the function $\hat J$ belongs to $\widehat {\cal W}_b$. Since by Prop.\ \ref{prp-proptaffmonft}(a) (or Prop.\ \ref{prp-proptaffmonfth}, depending on the assumptions), $\hat J$ is the unique solution of Bellman's equation within $\widehat {\cal W}$, it follows that $\hat J$ is the unique solution of Bellman's equation within $\widehat {\cal W}_b$. 

The proof of part (b) and that $\hat J$ is the unique solution of the Bellman Eq.\ \eqref{eq-bellmaneq} within the set $\widehat {\cal W}_b$ follow as in the proof of Prop.\ \ref{prp-proptaffmonft}. 

Assume now that $\hat J=\jstar$. Then from the preceding proof, $\jstar$ is the unique solution of Bellman's equation within the set $\widehat {\cal W}_b=\{J\in {\cal B}\mid \jstar\le J\}$. If there were another solution $J'$ within ${\cal B}$, then by Prop.\ \ref{prp-propnegdp}(a), we would have $\jstar\le J'$ so that $J'\in \widehat {\cal W}_b$. This shows that $J'=\jstar$, so $\jstar$ is the unique solution of Bellman's equation within ${\cal B}$.  \qed

\vspace{1pc}

The uniqueness of solution of Bellman's equation within  ${\cal B}$ when $\hat J=\jstar$ [cf.\ part (a) of the preceding proposition] is consistent with Example \ref{ex-example2}. In that example, $\jstar$ and $\hat J$ are equal and bounded, and all the additional solutions of Bellman's equation are unbounded. 

Note that without the assumption of existence of a uniformly proper $\p$, $\hat J$ and $\jstar$ need not belong to ${\cal B}$. As an example, let $X$ be the set of nonnegative integers, let $t=0$, and let there be a single policy that moves the system deterministically from a state $x\ge 1$ to the state $x-1$ at cost $g(x,x-1)=1$. Then 
$$\hat J(x)=\jstar(x)=x,\qquad \forall\ x\in X,$$
 so $\hat J$ and $\jstar$ do not belong to ${\cal B}$, even though $g$ is bounded. Here the unique policy is proper at all $x$, but is not uniformly proper. 
\vskip-1pc

\section{Concluding Remarks}

We have considered nonnegative cost SSP problems, which involve arbitrary state and control spaces, and a Bellman equation with possibly multiple solutions. Within this context, we have generalized the notion of a proper policy and we have discussed the restricted optimization over just the proper policies. 
The restricted optimal cost function $\hat J$ is a solution of Bellman's equation, and if the cost per stage is bounded, $\hat J$ is the maximal solution within the set of  nonnegative functions that are bounded within their effective domain.  By contrast, $\jstar$ is the minimal solution. When compared with their deterministic counterparts of the paper [3], the results of the present paper highlight an interesting difference: in deterministic problems $\hat J$ is the maximal solution of Bellman's equation within all functions in ${\cal J}$ (unbounded as well as extended real-valued), whereas this need not be true for stochastic problems.



\section{References}
\def\ref{\vskip1.pt\pn}

\def\refer{\ref}

\ref [1]  Bertsekas, D.\ P., and Shreve, S.\ E., 1978.\  Stochastic Optimal
Control:  The Discrete Time Case, Academic Press, N.\ Y.\ (republished by Athena Scientific, Belmont, MA, 1996); may be downloaded 
from\hfill\break http://web.mit.edu/dimitrib/www/home.html.

\ref[2] James, H.\ W., and Collins, E.\ J., 2006.\ ``An Analysis of Transient Markov Decision Processes," J.\
Appl.\ Prob., Vol.\ 43, pp.\ 603-621.

\ref[3] Bertsekas, D.\ P., 2018.\ ``Stable Optimal Control and Semicontractive Dynamic Programming," SIAM J.\ on Control and Optimization, Vol.\ 56, pp.\ 231-252.

\ref[4] Pallu de la Barriere, R., 1967.\ Optimal Control Theory, Saunders, Phila; republished by Dover, N. Y., 1980.

\ref [5] Derman, C., 1970.\ Finite State Markovian Decision Processes, Academic Press, N.\ Y.

\ref [6] Whittle, P., 1982.\  Optimization Over Time, Wiley, N.\ Y., Vol.\
1, 1982, Vol.\ 2, 1983.

\ref [7]  Bertsekas, D.\ P., and Tsitsiklis, J.\ N., 1989.\ Parallel and Distributed Computation: Numerical Methods, Prentice-Hall, Englewood Cliffs, N.\ Y.\ (republished by Athena Scientific, Belmont, MA, 1996); may be downloaded 
from\hfill\break http://web.mit.edu/dimitrib/www/home.html.

\ref [8] Puterman, M.\ L., 1994.\  Markov Decision Processes: Discrete Stochastic Dynamic Programming, J.\ Wiley, N.\ Y.

\ref[9] Altman, E., 1999.\ Constrained Markov Decision Processes, CRC Press, Boca Raton, FL.

\ref[10] Hernandez-Lerma, O., and Lasserre, J.\ B., 1999.\ Further Topics on Discrete-Time Markov Control Processes, Springer, N.\ Y.

\ref[11] Bertsekas, D.\ P., 2012.\ Dynamic Programming and Optimal Control, Vol.\ II: Approximate Dynamic Programming, Athena Scientific, Belmont, MA.

\ref [12]  Bertsekas, D.\ P., and Tsitsiklis, J.\ N., 1991.\ ``An Analysis of
Stochastic Shortest Path Problems,"
Math.\ of OR, Vol.\ 16, pp.\ 580-595.

\ref[13] Bertsekas, D.\ P., 2018.\ Abstract Dynamic Programming, 2nd Edition, Athena Scientific, Belmont, MA.

\ref[14] Bertsekas, D.\ P., 2015.\ ``Regular Policies in Abstract Dynamic Programming," Lab.\ for Information and Decision Systems Report LIDS-P-3173, MIT, May 2015; arXiv preprint arXiv:1609.03115; SIAM J.\ on Control and Optimization, Vol.\ 27, 2017, pp.\ 1694-1727.

\ref[15] Pliska, S.\ R., 1978.\ ``On the Transient Case for Markov Decision Chains with General State Spaces,"
in Dynamic Programming and its Applications, by M.\ L.\ Puterman
(ed.), Academic Press, N.\ Y.

\ref[16] Bertsekas, D.\ P., and Yu, H., 2016.\ ``Stochastic Shortest Path Problems Under Weak Conditions,"  Lab.\ for Information and Decision Systems Report LIDS-2909. 

\ref [17] Strauch, R., 1966.\  ``Negative Dynamic Programming," Ann.\ Math.\
Statist., Vol.\ 37, pp.\ 871-890.

\ref[18] Willems, J., 1971.\ ``Least Squares Stationary Optimal Control and the Algebraic Riccati Equation," IEEE Trans.\ on Automatic Control, Vol.\ 16, pp.\ 621-634.

\ref[19] Kucera, V., 1972.\ ``The Discrete Riccati Equation of Optimal Control," Kybernetika, Vol.\ 8, pp.\ 430-447.

\ref[20] Kucera, V., 1973.\ ``A Review of the Matrix Riccati Equation," Kybernetika, Vol.\ 9, pp.\ 42-61.

\vfill
 
\end{document}